\documentclass{amsart}

\usepackage{amsmath,amsthm,amssymb,amsfonts}
\usepackage{hyperref}

\input xy
\xyoption{all}
\usepackage{enumerate}

\makeatletter
\newcommand{\labitem}[2]{%
\def\@itemlabel{\textbf{#1}}
\item
\def\@currentlabel{#1}\label{#2}}
\makeatother

\theoremstyle{plain}
\newtheorem{thm}{Theorem}
\newtheorem{prp}[thm]{Proposition}

\theoremstyle{definition}

\theoremstyle{remark}
\newtheorem*{rmk}{Remark}

\newtheorem*{ack}{Acknowledgements}

\newcommand{\into}{\hookrightarrow}

\newcommand{\la}{\langle}
\newcommand{\ra}{\rangle}

\makeatletter
\@tempcnta=\@ne
\def\@nameedef#1{\expandafter\edef\csname #1\endcsname}
\loop\ifnum\@tempcnta<27
  \@nameedef{C\@Alph\@tempcnta}{\noexpand\mathcal{\@Alph\@tempcnta}}
  \advance\@tempcnta\@ne
\repeat

\makeatletter
\@tempcnta=\@ne
\def\@nameedef#1{\expandafter\edef\csname #1\endcsname}
\loop\ifnum\@tempcnta<27
  \@nameedef{B\@Alph\@tempcnta}{\noexpand\mathbb{\@Alph\@tempcnta}}
  \advance\@tempcnta\@ne
\repeat

\makeatletter
\@tempcnta=\@ne
\def\@nameedef#1{\expandafter\edef\csname #1\endcsname}
\loop\ifnum\@tempcnta<27
  \@nameedef{F\@Alph\@tempcnta}{\noexpand\mathsf{\@Alph\@tempcnta}}
  \advance\@tempcnta\@ne
\repeat

\begin{document}
\title[Conjugacy classes of cyclic subgroups]{Counting conjugacy classes of cyclic subgroups for fusion systems}

\author{Sejong Park}
\address{Section de math\'ematiques, \'Ecole Polythechnique F\'ed\'erale de Lausanne, Station 8, CH-1015 Lausanne, Switzerland}
\email{sejong.park@epfl.ch}

\date{\today}
\begin{abstract}
We give another proof of an observation of Th\'evenaz~\cite{T1989} and present a fusion system version of it.  Namely, for a saturated fusion system $\CF$ on a finite $p$-group $S$, we show that the number of the $\CF$-conjugacy classes of cyclic subgroups of $S$ is equal to the rank of certain square matrices of numbers of orbits, coming from characteristic bisets, the characteristic idempotent and finite groups realizing the fusion system $\CF$ as in our previous work~\cite{P2010}.
\end{abstract}
\maketitle

\section{Statements of the results} \label{S:intro}

%

In~\cite{T1989}, Th\'evenaz observed the `curiosity' that a finite cyclic group $G$ can be characterized by the nonsingularity of the matrix of the numbers of double cosets in $G$.  In fact He proved a more general fact that for an arbitrary finite group $G$ the number of the conjugacy classes of cyclic subgroups of $G$ is equal to the rank of that matrix.  This can be stated slightly more generally by introducing a subgroup $H$ of $G$ and considering the $G$-conjugacy classes of subgroups of $H$ as follows.

\begin{thm} \label{T:group}
Let $G$ be a finite group and let $H\leq G$.   The rank of the matrix
\[
	( |P\backslash G /Q|  )_{P,Q\leq_G H},
\]
whose rows and columns are indexed by the $G$-conjugacy classes of subgroups of $H$ and whose entries are the numbers of the corresponding double cosets in $G$, is equal to the number of the $G$-conjugacy classes of cyclic subgroups of $H$.
\end{thm}

In~\cite{P2010}, we observed that every saturated fusion system $\CF$ on a finite $p$-group $S$ can be realized by a finite group $G$ containing $S$ as a (not necessarily Sylow) $p$-subgroup.  Thus the above theorem yields a fusion system version as follows.

\begin{thm} \label{T:fusion-via-group}
Let $\CF$ be a saturated fusion system on a finite $p$-group $S$. Let $G$ be a finite group which contains $S$ as a subgroup and realizes $\CF$.  Then the rank of the matrix
\[
	( |P\backslash G /Q|  )_{P,Q\leq_G S}
\]
is equal to the number of the $\CF$-conjugacy classes of cyclic subgroups of $S$.
\end{thm}

By a result of Broto, Levi and Oliver \cite[Proposition 5.5]{BLO2003theory}, every saturated fusion system $\CF$ on a finite $p$-group $S$ has a (non-unique) characteristic biset $\Omega$.  See Section~\ref{S:fusion} for a precise definition; in particular, $\Omega$ is a finite $(S,S)$-biset, i.e., a finite set with compatible left and right $S$-actions.  If $\CF$ is the fusion system of  a finite group $G$ on its Sylow $p$-subgroup $S$, then $G$ is a characteristic biset for $\CF$ with the obvious $S$-action on the left and right.  So we may well expect that the matrix of the above theorem with $G$ replaced by $\Omega$ has the same rank.  Indeed this is the case.

\begin{thm} \label{T:fusion-via-biset}
Let $\CF$ be a saturated fusion system on a finite $p$-group $S$. Let $\Omega$ be a characteristic biset for $\CF$.  Then the rank of the matrix
\[
	( |P\backslash \Omega /Q|  )_{P,Q\leq_\CF S}
\]
of the number of $(P,Q)$-orbits of $\Omega$ indexed by the $\CF$-conjugacy classes of subgroups of $S$ is equal to the number of the $\CF$-conjugacy classes of cyclic subgroups of $S$.
\end{thm}

Finally, one can replace the characteristic biset $\Omega$ in the above theorem by the characteristic idempotent $\omega_\CF$ (which is a virtual $(S,S)$-biset; see Section~\ref{S:fusion}) with $|P\backslash \omega_\CF /Q|$ as the linearized number of $(P,Q)$-orbits. 

\begin{thm} \label{T:fusion-via-idempotent}
Let $\CF$ be a saturated fusion system on a finite $p$-group $S$. Let $\omega_\CF$ be the characteristic idempotent for $\CF$.  Then the rank of the matrix
\[
	( |P\backslash \omega_\CF /Q|  )_{P,Q\leq_\CF S},
\]
is equal to the number of the $\CF$-conjugacy classes of cyclic subgroups of $S$.
\end{thm}

We will give a proof of Theorem~\ref{T:group} (and hence obtain Theorem~\ref{T:fusion-via-group} as a corollary), which is slightly different from that of \cite{T1989}.  This new proof uses (at least explicitly) only the Burnside ring $B(G)$ of $G$, not the rational representation ring $R_\BQ(G)$ as in \cite{T1989}.  Therefore it is better suited for adapting to the fusion system case (Theorem~\ref{T:fusion-via-biset} and~\ref{T:fusion-via-idempotent}), which we do subsequently.


\section{The group case} \label{S:group}

We prove Theorem~\ref{T:group}.  As remarked  in Section~\ref{S:intro}, Theorem~\ref{T:fusion-via-group} then immediately follows as a corollary.

Let $G$ be a finite group.  Let $B(G)$ be the Burnside ring of $G$, i.e., the Grothendieck ring of the isomorphism classes $[X]$ of finite $G$-sets $X$.  As an additive group, $B(G)$ is a free abelian group with the canonical basis $\{ [G/P] \mid P\leq_G G \}$.  Let $\BQ B(G) = \BQ \otimes_\BZ B(G)$ and regard $B(G)$ as a subgring of $\BQ B(G)$. In particular the canonical basis for $B(G)$ is a $\BQ$-basis for $\BQ B(G)$.

It is a well-known fact that for each $P\leq G$ the fixed-point map
\[
	\chi_P\colon B(G) \to \BZ,\quad [X] \mapsto |X^P|,
\]
is a ring homomorphism which depends only on the $G$-conjugacy class of $P$, and their product (tensored with $\BQ$) 
\[
	\chi= \BQ \otimes_\BZ \prod_{P\leq_G G} \chi_P \colon \BQ B(G) \to \prod_{P\leq_G G} \BQ
\]
is a $\BQ$-algebra isomorphism. For each $P\leq G$, let $e^G_P$ denote the element of $\BQ B(G)$ such that
\[
	\chi_Q(e^G_P) = 
	\begin{cases}
		1,	&\text{$P=_G Q$},\\
		0,	&\text{otherwise}.
	\end{cases}
\]
Then again the element $e^G_P$ depends only on the $G$-conjugacy class of $P$ and $\{ e^G_P \mid P \leq_G G \}$ is a set of pairwise orthogonal primitive idempotents of $\BQ B(G)$ whose sum is equal to $1$; in particular it is a $\BQ$-basis for $\BQ B(G)$.  Furthermore, for $H\leq G$, let $B(G)_H$ be the subgroup of $B(G)$ generated by the elements $[G/P]$ with $P\leq_G H$.  Then $\BQ B(G)_H = \BQ \otimes_\BZ B(G)_H$ is a subalgebra of $\BQ B(G)$ with $\BQ$-basis $\{ [G/P] \mid P \leq_G H \}$.  Note that the elements $e^G_P$ with $P\leq_G H$ belong to $\BQ B(G)_H$ and hence  $\{ e^G_P \mid P \leq_G H \}$ is another basis for $\BQ B(G)_H$.

For each $P\leq G$ consider the $\BQ$-linear map
\[
	\rho_P \colon \BQ B(G) \to \BQ,\quad [X] \mapsto |P\backslash X|,
\]
which counts the $P$-orbits.  By Burnside's orbit counting lemma, we have
\[
	\rho_P(x) = \frac{1}{|P|} \sum_{u\in P} \chi_{\la u \ra}(x),\quad x\in\BQ B(G).
\]	
Thus
\begin{equation}\label{E:orbit-counting}
	\rho_P(e^G_Q) \neq 0 \iff \text{$Q$ is cyclic and $Q\leq_G P$}.
\end{equation}
Now the given matrix in Theorem~\ref{T:group} is equal to
\[
	(\rho_P(G/Q))_{P,Q\leq_G H}. 
\]
By change of basis, this matrix has the same rank as
\[
	(\rho_P(e^G_Q))_{P,Q\leq_G H}. 
\]
List the subgroups of $H$ (up to $G$-conjugacy) in two blocks, the first consisting of cyclic subgroups and the second of noncyclic subgroups, and with nondecreasing order in each block.  Then by~\eqref{E:orbit-counting} the above matrix has the form
\[
	\begin{pmatrix}
		A & 0 \\
		B & 0
	\end{pmatrix}
\]
where $A$ is a lower triangular matrix with nonzero diagonal entries. Thus Theorem~\ref{T:group} follows.

\section{The fusion system case} \label{S:fusion}

We first prove Theorem~\ref{T:fusion-via-biset}.  In fact, we prove a slightly generalized version of it.

\begin{prp} \label{T:fusion-via-biset-general}
Let $\CF$ be a saturated fusion system on a finite $p$-group $S$. Suppose that $\Omega$ is a finite $(S,S)$-biset which is $\CF$-stable and $\CF$-generated and which contains the obvious $(S,S)$-biset $S$. Then the rank of the matrix
\[
	( |P\backslash \Omega /Q|  )_{P,Q\leq_\CF S}
\]
is equal to the number of the $\CF$-conjugacy classes of cyclic subgroups of $S$.
\end{prp}

We first explain the terminology.  Let $\CF$ be a saturated fusion system on a finite $p$-group $S$.  An $S$-set $X$ is {\em $\CF$-stable} if, for all $P\leq S$ and all $\CF$-morphism $\varphi\colon P\to S$, the restrictions of the $S$-action on $X$ to $P$ via the inclusion $P\into S$ and via $\varphi\colon P\to S$ give isomorphic $P$-sets.  We say that an $(S,S)$-biset is {\em $\CF$-stable} if it is $\CF\times\CF$-stable viewed as a left $S\times S$-set by inverting the right action of $S$.  An $(S,S)$-biset is {\em $\CF$-generated} if, viewed as a left $S\times S$-set, all its isotropy subgroups are of the form $\Delta(P,\varphi) = \{ (u,\varphi(u)) \mid u\in P \}$ with $P\leq S$, $\varphi\colon P\to S$ in $\CF$.  A finite $(S,S)$-biset $\Omega$ is called a {\em characteristic biset} for $\CF$ if it is $\CF$-stable and $\CF$-generated and such that $|\Omega|/|S|$ is not divisible by $p$.  It is easy to see that every characteristic biset $\Omega$ contains the $(S,S)$-biset $S$.

Define
\[
	B(\CF) = \{ x \in B(S) \mid \chi_{P}(x) = \chi_{P'}(x) \text{ for all $P, P'\leq S$ with $P=_\CF P'$} \}.
\]	
Clearly $B(\CF)$ is a subring of $B(S)$, which is called the Burnside ring of the fusion system $\CF$.  For a finite $S$-set $X$, we have $[X] \in B(\CF)$ if and only if $X$ is $\CF$-stable. %
As before let $\BQ B(\CF) = \BQ \otimes_\BZ B(\CF)$.  Clearly the elements
\[
	e^\CF_P := \sum_{P'=_\CF P} e^S_{P'},
\]
where $P\leq S$ and the sum is over the $S$-conjugacy classes of subgroups $P'$ of $S$ which are $\CF$-conjugate to $P$, belong to $\BQ B(\CF)$.  The set $\{ e^\CF_P \mid P\leq_\CF S \}$ is a set of pairwise orthogonal primitive idempotents of $\BQ B(\CF)$ whose sum is equal to $1$; in particular it is a $\BQ$-basis for $\BQ B(\CF)$.  By \eqref{E:orbit-counting}, we have
\begin{equation}\label{E:orbit-counting-F}
	\rho_P(e^\CF_Q) \neq 0 \iff \text{$Q$ is cyclic and $Q\leq_\CF P$}.
\end{equation}

Let $\Omega$ be the $(S,S)$-biset given in the above proposition.  By the $\CF$-stability of $\Omega$, the left $S$-set $\Omega/P$ of the right $P$-orbits of $\Omega$ is also $\CF$-stable for $P\leq S$.  Moreover 
\[
	\chi_Q([\Omega/P]) \neq 0 \implies Q \leq_\CF P;\quad \chi_P([\Omega/P]) \geq |N_S(P)/P|.
\]
The former follows from that $\Omega$ is $\CF$-generated and the latter from that $\Omega$ contains $S$. Hence
\[
	\{ [\Omega/P] \mid P\leq_\CF S \}
\]
is a $\BQ$-basis for $\BQ B(\CF)$. Thus the matrix
\[
	( |P\backslash \Omega /Q|  )_{P,Q\leq_\CF S} = ( \rho_P([\Omega /Q])  )_{P,Q\leq_\CF S}
\]
has the same rank as 
\[
	( \rho_P(e^\CF_Q)  )_{P,Q\leq_\CF S},
\]
which is equal to the number of the $\CF$-conjugacy classes of cyclic subgroups of $S$ by \eqref{E:orbit-counting-F}.

\begin{rmk}
Note that the finite group $G$ in Theorem~\ref{T:fusion-via-group}, viewed as an $(S,S)$-biset, satisfies the hypotheses for $\Omega$ in Proposition~\ref{T:fusion-via-biset-general}.  Thus Theorem~\ref{T:fusion-via-group} can also be obtained from Proposition~\ref{T:fusion-via-biset-general}.
\end{rmk}

Now we address Theorem~\ref{T:fusion-via-idempotent}.  In Proposition~\ref{T:fusion-via-biset-general}, the condition that $\Omega$ contains the $(S,S)$-biset $S$ is equivalent to that $\chi_P(\Omega/P)\neq 0$ for all $P\leq S$, given the other conditions on $\Omega$. Proposition~\ref{T:fusion-via-biset-general} then applies to all virtual $(S,S)$-bisets $\omega$ with coefficients in $\BQ$ which are $\CF$-stable, $\CF$-generated and such that $\chi_P(\omega/P) \neq 0$ for all $P\leq S$, where $\omega/P$ denotes the linearized right $P$-orbits of $\omega$.  The proof is identical to the one given above.  In particular, Reeh~\cite[Proposition 4.5, Corollary 5.8]{Reeh2013Transfer} shows that if $\omega$ is the {\em characteristic idempotent} of $\CF$, i.e., the unique virtual $(S,S)$-biset with coefficients in $\BZ_{(p)}$ which is $\CF$-stable, $\CF$-generated and which is an idempotent in the double Burnside ring $\BZ_{(p)}B(S,S)$, then the elements $\omega/P = \omega \circ_S [S/P] = \beta_P$ with $P\leq_\CF S$ form a basis of $\BZ_{(p)}B(\CF)$ such that $\chi_P(\omega/P) \neq 0$.  This proves Theorem~\ref{T:fusion-via-idempotent}.

\begin{ack}
We thank Prof.\ Markus Linckelmann for pointing us to the  paper of Th\'evenaz and Prof.\ Th\'evenaz for helpful discussions.
\end{ack}



\begin{thebibliography}{1}

\bibitem{BLO2003theory}
C.~Broto, R.~Levi and B.~Oliver, \emph{The homotopy theory of fusion systems}, J. Amer. Math. Soc. \textbf{16} (2003), no.~4, 779--856
  (electronic).

\bibitem{P2010} S.~Park, \emph{Realizing a fusion system by a single finite group}, Arch. Math. (Basel) 94 (2010), no. 5, 405-410.

\bibitem{Reeh2013Transfer}
S.~Reeh, \emph{Transfer and characteristic idempotents for saturated fusion systems}, preprint (2013).

\bibitem{T1989}
J.~Th{\'e}venaz, \emph{A characterization of cyclic groups}, Arch. Math. (Basel) \textbf{52} (1989), no.~3, 209--211.

\end{thebibliography}
\end{document}